\documentclass[12pt,a4paper,onecolumn,twoside]{article}

\usepackage{geometry}
\usepackage[english]{babel}
\usepackage[latin1]{inputenc}
\usepackage[T1]{fontenc}
\usepackage{fancyhdr}        

\usepackage{amsmath}
\usepackage{amsfonts}
\usepackage{amssymb}
\usepackage{amsthm}          
\usepackage{graphicx}        

\usepackage{verbatim}
\usepackage[all]{xy}         
\usepackage{syntonly}
\usepackage[pdftex,bookmarks,colorlinks,citecolor=blue,urlcolor=blue]{hyperref}          
\usepackage{stackrel}        

\newcommand\NN{\mathbb{N}}   
\newcommand\ZZ{\mathbb{Z}}   
\newcommand\II{\mathbb{I}}   
\newcommand\RR{\mathbb{R}}   
\newcommand\CC{\mathbb{C}}   

\newcommand\TF{\mathcal{F}}     
\newcommand{\HH}{\mathcal{H}}   
\newcommand{\BB}{\mathcal{B}}   
\newcommand{\VV}{\mathcal{V}}   

\newcommand{\rr}{\mathcal{R}}   

\newcommand{\VS}        {vector space}
\newcommand{\NS}        {normed space}
\newcommand{\BS}        {Banach space}
\newcommand{\BSS}       {Banach spaces}
\newcommand{\HS}        {Hilbert space}

\newcommand{\RKHS}      {reproducing kernel Hilbert space}
\newcommand{\RKHSS}     {reproducing kernel Hilbert spaces}
\newcommand{\RKBS}      {reproducing kernel Banach space}
\newcommand{\RKBSS}     {reproducing kernel Banach spaces}
\newcommand{\UDF}       {uniformly Fr\'echet differentiable}

\newcommand{\SIPRKBS}   {s.i.p.{ }RKBS}
\newcommand{\UC}        {uniformly convex}

\newcommand{\UU}        {uniform}

\newcommand{\PI}        {inner product}
\newcommand{\SIP}       {semi-inner product}
\newcommand{\RK}        {reproducing kernel}
\newcommand{\SIPK}      {s.i.p. kernel}
\newcommand{\SIPRK}     {s.i.p. reproducing kernel}

\newcommand{\BON}       {orthonormal basis}
\newcommand{\SB}        {Schauder basis}
\newcommand{\RB}        {Riesz basis}
\newcommand{\BIO}       {biorthogonal}

\newcommand{\KST}       {Kramer sampling theorem}
\newcommand{\RRT}       {Riesz representation theorem}

\newcommand{\IFF}       {if and only if}

\newtheorem{theo}{Theorem}[section]
\newtheorem{prop}[theo]{Proposition}

\newtheorem{coro}[theo]{Corollary}

\theoremstyle{definition}
\newtheorem{deff}[theo]{Definition}
\newtheorem*{exam}{Example}

\theoremstyle{remark}
\newtheorem{comm}[theo]{Remark}

\title{
Sampling basis in \RKBSS}
\date{}
\begin{document}

\lhead[\thepage]{\nouppercase{\leftmark}}                                   
\chead[]{}
\rhead[\nouppercase{\rightmark}]{\thepage}                                    
\renewcommand{\headrulewidth}{0pt}

\lfoot[]{}
\cfoot[]{}
\rfoot[]{}
\renewcommand{\footrulewidth}{0pt}

\fancypagestyle{plain}{
\fancyhead[L]{}
\fancyhead[C]{}
\fancyhead[R]{}
\fancyfoot[L]{}
\fancyfoot[C]{\thepage}
\fancyfoot[R]{}
\renewcommand{\headrulewidth}{0pt}
\renewcommand{\footrulewidth}{0pt}     
}

\pagestyle{fancy} 
\author{{ Hern\'{a}n Centeno}\\ 
{ {\normalsize   Universidad de Buenos Aires,}}\\ 
{{\normalsize Facultad de Ciencias Exactas y Naturales}}\\ 
{{\normalsize y Facultad de Ingenier\'{\i}a}}\\ 
{{\normalsize Departamento de Matem\'atica}}\\
{{\normalsize  Paseo Col\'on 850, CABA, CP. 1063 , Argentina}}\\
{ {\normalsize  hcenteno@dm.uba.ar }}\\
\\
{ {\normalsize Juan Miguel Medina}}\\
{ {\normalsize  Universidad de Buenos Aires,}}\\{{ Departamento de Matem\'atica}}\\ 
{{Facultad de Ingenier\'{\i}a y IAM-CONICET}}\\
{ {\normalsize  Paseo Col\'on 850, CABA, CP. 1063 , Argentina}}\\
{ {\normalsize  jmedina@fi.uba.ar }}}

\maketitle
%
\begin{abstract}
We present necessary and sufficient conditions to hold true a Kramer type  sampling theorem  over \SIP\ \RKBSS. Under some sampling-type hypotheses over a sequence of functions on these \BSS\  
it results necessary that such sequence must be a $X_d$-\RB\ and a sampling basis for the space. These results are a generalization of some already known sampling theorems over \RKHSS. 
\end{abstract}
%
%
\begin{flushleft}
\end{flushleft}
%
%
\section{Introduction}
\markboth{\thesection. Introduction}{\thesection. Introduction}
%
The celebrated sampling theorem of Whittaker-Shannon-Kotel'nikov (1933) \cite{Eldar, Zayed} establishes that all
finite energy function $f\in L^2 (\RR)$ band-limited to $\left[ -\sigma,\sigma \right]$, i.e., 
the \textit{Fourier transform} of $f$ is supported on the interval $\left[ -\sigma,\sigma \right]$, 
can be completely recovered through samples in the integers $\{f(n)\}_{n\in\ZZ}$, obtaining in this way the following representation
\begin{equation}  
                f(t)= \sum_{n=-\infty}^{\infty} f\bigg(\frac{n}{2\sigma} \bigg ) \dfrac{\sin \pi( 2\sigma t-n)}{\pi (2\sigma t-n)}
                \qquad t\in\RR \nonumber
\end{equation}
with the series being absolutely and uniformly convergent on compact subsets of $\RR$. By writing it a bit different, 
we note that the band-limited functions can be given by
\begin{equation}  
                f(t)= \frac{1}{\sqrt{2\pi} }\int_{-\sigma}^{\sigma} F(x)e^{- it\omega} d\omega
                \qquad t\in \RR  ,\nonumber
\end{equation}
being $\{e^{ i n(\cdot )}\}_{n\in\ZZ}$ an \BON\ of $L^2 \left[ -\sigma,\sigma \right] $.
By noting this, later in the 1959, Kramer \cite{Eldar, Kram} extended this result to functions defined by another integral operator 
$TF=f$, now with kernel $\kappa$ instead of the exponentials:
\begin{displaymath}
                   f(t)=\int_I F(x)\overline{\kappa (t,x)} dx
                   \qquad t\in \RR
\end{displaymath}
where $I$ is a compact interval of $\RR$ and $\kappa (t,\cdot)\in L^2(I)\;\;\forall\,t\in \RR$. The existence of a sequence 
$\{t_n\}_{n\in\ZZ}\subset \RR$ such that $\{ \kappa (t_n,\cdot) \}_{n\in\ZZ}$ is an orthogonal and complete sequence in $L^2(I)$ 
was the hypothesis used by Kramer for this result to hold. Thanks to this he obtained the sampling expansion for such functions:
\begin{equation}
                f(t)= \sum_{n=-\infty}^{\infty} f(t_n ) S_{n}(t)
                      \quad \textrm{with} \quad
                S_{n}(t)= \frac{\int_{I}\kappa(t_n,x)\overline{\kappa(t,x)} dx}{\int_{I}\vert \kappa(t_n,x)\vert^2 dx}
                      \qquad t\in\mathbb{R} \nonumber
\end{equation}
as before, the series is absolutely convergent. This result allow us to work in \textit{non uniform} sampling problems in contrast to 
the Whittaker-Shannon-Kotel'nikov sampling theorem.
Both of the integral operators could be written by using the usual inner product of $L^2(I)$ and then we obtain a 
possible direction to where it can be generalized this \KST.

Thanks to the theory of \textit{\RKHSS} (written RKHS for short) by Aronszajn \cite{Aro} in the 1950 and its particular case of functions 
which are image by an integral operator (Saitoh 1988, \cite{Sai}), the previous sampling results can be naturally viewed inside this 
framework. Thanks to a new generalization (again by Saitoh), it can be considered like particular cases of the so-called 
\textit{Abstract \KST} (Garc\'ia, Hern\'andez-Medina \& Mu\~{n}oz-Bouzo, 2014 \cite{Gar3}), where the functions now have the form:
\begin{equation}  
                f(t) = \langle x , \Phi (t) \rangle
                \qquad t\in \Omega  ,\nonumber
\end{equation}
where $\Omega$ is an arbitrary set, $(\HH , \langle\cdot , \cdot\rangle )$ is a \textit{\HS} and $\Phi:\Omega\rightarrow \HH$ is an 
arbitrary function. Under the hypotheses of the existence of sequences 
$\{t_n\}_{n\in\NN}\subset \Omega$, $\{a_n\}_{n\in\NN}\subset\CC\setminus\{0\}$ and $\{x_n\}_{n\in\NN}\subset \HH$ a \textit{\RB} such that 
the sequence $\{ \Phi (t_n) \}_{n\in\NN}$ satisfies the interpolation condition $\Phi (t_n)= \overline{a_n} x_n \;\;\forall\,n\in\NN$, 
they were able to prove that
\begin{equation}
                f(t)= \sum_{n=1}^{\infty} f(t_n ) \frac{S_{n}(t)}{a_n}
                      \quad \textrm{with} \quad
                S_{n}(t)= \langle y_n , \Phi (t) \rangle
                      \qquad t\in\Omega \nonumber
\end{equation}
where $\{y_n\}_{n\in\NN}\subset \HH$ is the \BIO\ \RB\ of $\{x_n\}_{n\in\NN}$ and the series is convergent in the RKHS-norm that 
contains such functions, also, the convergence is absolute and uniform on subsets of $\Omega$ where the map 
$t\mapsto \Vert \Phi (t) \Vert$ is bounded. 

Due to the recent theory of \RKBSS\ (written RKBS for short) developed by Zhang, Xu \& Zhang \cite{Zhang1} and the subsequent theory of 
$X_d$-Bessel sequences, $X_d$-frames and $X_d$-\RB\ by Zhang \& Zhang \cite{Zhang2}, Garc\'ia \& Portal (2013, \cite{Gar2}) were able 
to extend the last result (stated in Section \ref{sec:krasampling}) to the \BSS\ setting. By using these recent concepts we state and 
prove a generalization of the following possible ``{\textit{converse}}'' of the \KST:

\begin{theo}[A converse of the \KST\ \cite{Gar1}]\label{conversekst}
Let $\HH$ be the range of the integral linear transform $T : L^2(I) \ni F \rightarrow f \in \HH$ considered as a RKHS with the kernel 
$k$ defined by $k(t, s) := \langle K(\cdot , t),K(\cdot , s)\rangle_{L^2(I)}$. Let $\{ S_n\}_{n=0}^{\infty}$ be a sequence in $\HH$ such 
that $\sum_{n=0}^{\infty} \vert S_n(t)\vert^2 <+\!\infty$, $t\in\Omega$ and let $\HH_{samp}$ be a RKHS corresponding to the kernel 
$K_{samp}(s, t) := \sum_{n=0}^{\infty}\overline{S_n(s)}S_n(t)$. Then, we have the following results:
\begin{itemize}
\item[$1^\circ$)] Suppose that the sequence $\{ S_n\}_{n=0}^{\infty}$ satisfies the condition that for each sequence 
$\{ \alpha_n\}_{n=0}^{\infty}\in\ell_2 (\NN_0)$ such that $\sum_{n=0}^{\infty}\alpha_n S_n(t)=0$ implies $\alpha_n = 0$ for all $n$. 
Then, $\HH_{samp}\subset\HH$ and $\{ S_n\}_{n=0}^{\infty}$ is an \BON\ in $\HH_{samp}$.

\item[$2^\circ$)] Suppose in addition to $1^\circ$) the existence of sequences $\{ t_n\}_{n=0}^{\infty}$ in $\Omega$ and 
$\{ a_n\}_{n=0}^{\infty}$ in $\CC\setminus \{0\}$ such that
\begin{equation}
                \bigg \{ \frac{f(t_n )}{a_n} \bigg \}_{n\in\NN_{0}} \in \ell_2 (\NN_{0})
                \quad \textnormal{ and } \quad
                f(t) =\sum_{n=0}^{\infty} f(t_n)\frac{S_n (t)}{a_n}
                \qquad \textnormal{ for any } f\in\HH \nonumber
\end{equation}
where the sampling series is pointwise convergent in $\Omega$. Then
\begin{itemize}
\item[$\bullet$] $\HH_{samp} = \HH$.
                                                                    
\item[$\bullet$] The norms of $\HH_{samp}$ and $\HH$ are equivalent, i.e., for some constants $0<a\leq b$
                 \begin{equation}
                                 a\Vert f \Vert_{samp} \leq \Vert f \Vert_{\HH} \leq b\Vert f \Vert_{samp} \nonumber
                 \end{equation}
                 Consequently $\{ S_n\}_{n=0}^{\infty}$ is a \RB\ for $\HH$.
                                                                          
\item[$\bullet$] The sequences $\{\overline{a_i^{-1}}K(\cdot ,t_i)\}_{i=0}^{\infty}$ and 
                 $\{ \sum_{n=0}^{\infty}\langle S_j,S_n\rangle_{\HH}K(\cdot ,t_n) \}_{j=0}^{\infty}$ as well as the sequences $\{ S_i\}_{i=0}^{\infty}$ and 
                 $\{ \sum_{n=0}^{\infty} k_{t_j}(t_n) a_n^{-1}S_n\}_{j=0}^{\infty}$ are biorthonormal in $L^2(I)$ and $\HH$ respectively.

\item[$\bullet$] If $a = b$ then $a^2 k(s, t) = k_{samp}(s, t)$ for all $s, t\in\Omega$ and the sequence $\{ S_n\}_{n=0}^{\infty}$ is a complete and orthogonal set 
                 in $L^2(I)$.                                  
\end{itemize}
\end{itemize}
\end{theo}

Recently, in \cite{Higg} is obtained another possible converse with different choices of hypotheses.
In the next section we give the preliminaries needed for the extension of this theorem to the \BS\ setting. We only list the results 
and invite to the reader to see \cite{dragomir,Giles,Lumer,Zhang1,Zhang2} for much more details.
%
\section{Definitions and basic results}
%
%
%
\subsection{The normalized duality mapping and \SIP s}
Let ($E , \Vert\cdot \Vert$) be a \NS\ over $\CC$ and ($E^* , \Vert\cdot \Vert_{*}$) its corresponding dual space formed by the 
$\Vert\cdot \Vert$-continuous $\CC$-linear functional. We have defined the bilinear form 
$(\cdot , \cdot )_{E} : E \times E^* \rightarrow \CC$ given by $(f,f^*)_E= f^* (f),\ f\in X$ and $f^*\in E^*$. The mapping 
$J:E\rightarrow 2^{E^*}$ given by
\begin{equation}\label{eq:mapa-dual-def}
                J(f) = \{ f^*\in E^* : f^* (f)= \Vert f \Vert \Vert f^* \Vert_{*},\ \Vert f \Vert =\Vert f^* \Vert_{*} \}
                \qquad f\in E  \nonumber
\end{equation}
will be called \textit{the normalized duality mapping of the \NS} $E$ or shortly \textit{the dual map of} $E$. For our purposes, 
here and henceforth $E$ will be a \UU\ \BS, i.e., \UDF\ and \UC\ space \cite{Megg}. In this way, given $f\in E$ there exists a unique 
$f^*\in E$ such that $J(f)=\{f^*\}$ and so we have an isometric bijection $f\mapsto f^*$ between $E$ and $E^*$. For the proofs of 
these statements and more about the dual map see for example \cite{dragomir} and the references therein. We introduce the \SIP s (s.i.p. 
for short), these share almost all properties of the \PI s.

\begin{deff}\label{def-sipspace} 
Let $\VV$ be a $\CC$-\VS, a map $[\cdot\, , \cdot ] :\VV \times \VV \rightarrow \CC$ is called a \textit{\SIP} 
(in Lumer's sense \cite{Lumer}) if $\forall\, \alpha\in\CC$ and $\forall\, x,y,z\in \VV$ satisfies:

$ \bullet\quad [\alpha x+y,z ] = \alpha [ x,z ] + [ y,z ]$. 

$\bullet\quad [ x,\alpha  y ] = \overline{\alpha} [ x,y ]$.    
                                       
$\bullet\quad [ x,x ] > 0$, if $x\neq 0$. 

$\bullet\quad \vert [ x,y ] \vert^{2} \leq [ x,x ] [ y,y ]$. 
\end{deff}

When $E$ is a \UU\ \BS\ (so it is $E^*$) there exists a unique s.i.p. $[\cdot\, , \cdot ]$ on $E$ (hence a unique s.i.p. 
$[\cdot\, , \cdot ]_*$ on $E^*$) which is compatible with the norm in the sense that $\Vert f \Vert^2=[f,f]\ \, \forall\, f\in E$, also 
we have a \RRT, concretely, for each $L\in E^*$, there exists a unique $f\in E$ such that $L=f^*$ and $f^*(g)=[g,f]\ \,\forall \,g\in E$. 
The relationship between both \SIP s is given by $[ f^*,g^* ]_*=[ g,f ]\ \, f,g\in E$.
\subsection{Bessel sequences, Frames and Riesz bases via s.i.p.}
The following are included in \cite{Zhang2}. A \textit{BK-space} $X_d$ on a countable well-ordered index set $\II$ is a \BS\ of sequences 
indexed by $\II$ where the canonical vector forms a \SB. We impose the following additional conditions over $X_d$: it is a reflexive 
space, which guarantees its dual $X_d^*$ is also a BK-space and the duality between them is given by 
$( c, d )_{X_d} = \sum_{j\in\II} c_j d_j,\ \forall\, c = \{ c_j \}_{j\in\II} \in X_d ,\, d = \{ d_j \}_{j\in\II} \in X_{d}^{*}$; if the 
series $\sum_{j\in\II} c_j d_j$ converges in $\CC$ for all $c \in X_d$ then $d \in X_d^*$ and vice versa; finally the series 
$\sum_{j\in\II} c_j d_j$ converges absolutely in $\CC$ for all sequences $c \in X_d ,\, d  \in X_{d}^{*}$. 
For another types of sequence spaces we refer to \cite{Stoeva1, Stoeva2}. 

Given a sequence $\{f_{j}\}_{j\in\II}$ in $E$ we note by $\{f_{j}^*\}_{j\in\II}$ its dual sequence in $E^*$. A sequence 
$\{f_{j}\}_{j\in\II}$ in $E$ is called \textit{minimal}, if $f_{k}\notin\overline{span} \{ f_j : k\neq j \}\;\forall\, k\in\II$ and is 
called \textit{complete}, if $\overline{span} \{ f_{j} : j\in\II \} = E$. We have the following characterizations:

\begin{prop}\label{minimaliff}
Let $\{f_{j}\}_{j\in\II}$ be a sequence in $E$, then:

a) $\{f_{j}\}_{j\in\II}$ is minimal \IFF\ $\exists\, \{g_{j}\}_{j\in\II}$ in $E$ such that 
$[f_j , g_k]=\delta_{j,k}\ \forall\, j,k\in\II$.

b) $\{f_{j}\}_{j\in\II}$ is complete \IFF\ $f\in E$ is such that $[f_j ,f]=0\ \forall\, j\in\II$ then $f=0$.
\end{prop}
Where $\delta_{j,k}$ denotes the Kronecker's delta. The sequence $\{g_{j}\}_{j\in\II}$ in $a)$ is called a \BIO\ sequence of 
$\{f_{j}\}_{j\in\II}$ and when $\{f_{j}\}_{j\in\II}$ is also a complete sequence in $E$, then $\{g_{j}\}_{j\in\II}$ is unique.

We give first the definition of $X_d$-\textit{Riesz-Fischer} sequences and then introduce $X_d$-\textit{Bessel} sequences, 
$X_d$-\textit{frames} and $X_d$-\textit{\RB} at the same time as its characterizations, these will be used in the main 
result in Section \ref{sec:mainresult}. See \cite[Proposition 2.3 -- 2.13]{Zhang2}.

\begin{deff}[$X_d$-\textit{Riesz-Fischer sequences}]\label{rieszfischerseqiff}
$\{f_{j}\}_{j\in\II}\subset E$ is a $X_d$-\textit{Riesz-Fischer} sequence for $E$ if
\begin{equation}\label{eq:xd-riesz-fischer-def}  
\forall\, c=\{ c_j\}_{j\in\II} \in X_d ,
\quad \exists\,f\in E \quad \textnormal{ such that } \quad
[ f,f_j ]=c_j 
\quad\forall\, j\in\II .
\end{equation}
\end{deff}

\begin{prop}[$X_d$-\textit{Bessel} sequences]\label{besselseqiff}
Let $\{f_{j}\}_{j\in\II}$ be a sequence in $E$, are equivalent:

{\ \ }i) ($X_d$-\textit{Bessel} definition) There exists a constant $B>0$ such that
         \begin{equation}\label{eq:xd-bessel-def} 
         \Vert \{[ f,f_j ] \}_{j\in\II} \Vert_{X_d} \leq B\Vert f \Vert_{E} 
         \qquad \forall\, f\in E 
         \end{equation}

{\ }ii)  $U: E \rightarrow X_d$ given by $Uf=\{ [ f,f_j ] \}_{j\in \II}\ f\in E$ is a well-defined bounded operator.

iii)     $U^* :X_d^* \rightarrow E^*$ given by
         \begin{equation}\label{eq:synthesis}
         U^* d =\sum_{j\in\II} d_j f_j^* 
         \qquad\forall\, d=\{ d_j \}_{j\in\II}\in X_d^*
         \end{equation}

$\quad$ is a bounded operator and the series $\sum_{j\in\II} d_j f_j^*$ converges unconditionally in $E^*$.
\end{prop}

\begin{prop}[$X_d$-\textit{frames}]\label{framesiff}  
Let $\{f_{j}\}_{j\in\II}$ be a sequence in $E$, are equivalent:

{\ \ }i) ($X_d$-\textit{frame} definition) There exists constants $B\geq A>0$ such that
         \begin{equation}\label{eq:xd-frame-def} 
         A\Vert f \Vert_{E} \leq \Vert \{[ f,f_j ] \}_{j\in\II} \Vert_{X_d} \leq B\Vert f \Vert_{E} 
         \qquad \forall\, f\in E 
         \end{equation}

{\ }ii)  $U: E \rightarrow X_d$ is bounded and bounded below.

iii)     $U^* : X_d^* \rightarrow E^*$ is bounded and surjective.
\end{prop}
It is clear that an $X_d$-\textit{frame} for $E$ is an $X_d$-\textit{Bessel} sequence for $E$.

\begin{prop}[$X_d$-\textit{\RB}]\label{rieszbasisiff}
Let $\{f_{j}\}_{j\in\II}$ be a sequence in $E$, are equivalent:

{\ \ }i) ($X_d$-\textit{\RB} definition) $\{f_{j}\}_{j\in\II}$ is complete and $\exists\, B\geq A>0$ such that 
         \begin{equation}\label{eq:baseriesz-def}
         A\Vert c \Vert_{X_d} \leq \bigg \Vert \sum_{j\in\II} c_j f_j \bigg \Vert_{E} \leq B\Vert c \Vert_{X_d} 
         \qquad \forall\, c = \{ c_j \}_{j\in\II} \in X_d ,
         \end{equation}

{\ }ii)  $\{f_{j}^*\}_{j\in\II}$ is an $X_d^*$-\textit{frame} for $E^*$ and $\{f_{j}\}_{j\in\II}$ is a minimal sequence in $E$. 

iii)     $\{f_{j}\}_{j\in\II}$ is complete and $V : E^* \rightarrow X_d^*$ is bounded and surjective.

iv)      $\{f_{j}\}_{j\in\II}$ is complete and $V^* : X_d \rightarrow E$ is bounded and bounded below.
\end{prop}
\subsection{Reproducing kernel Banach spaces}

A \RKHS\ on a set $\Omega$ is a \HS\ $(\HH, \langle \cdot , \cdot\rangle )$ of $\CC$-valued functions on $\Omega$ and the point evaluations 
in $t\in\Omega$ are continuous linear functionals on $\HH$. The second condition is equivalent to the existence of a function 
$K:\Omega\times\Omega\rightarrow\CC$ such that $K(t,\cdot)\in\HH$ for each $t\in\Omega$, and for each $f\in\HH$ there holds the \textit{reproducing property}:
\begin{equation}\label{eq:propiedad-reproductiva-h}
                f(t)= \langle f,K(t,\cdot)\rangle
                \qquad t\in\Omega  .
\end{equation}
where the choice of the first variable of $K$ is simply by convenience (the second one is usually used). $K$ is unique and is called the \textit{\RK} for $\HH$. For our main purpose of doing sampling theory,  we adopt the next definition of \RKBS\ \cite{Zhang1} to extend 
these Hilbert spaces to the \BS\ setting.
                                                                     
\begin{deff}\label{rkbs}
A \textit{\RKBS} on a set $\Omega$ is a reflexive \BS\ $(\BB , \Vert\cdot\Vert )$ of $\CC$-valued functions on $\Omega$ for which
$\BB^*$ is isometrically isomorphic to a \BS\ $\BB^\#$ of $\CC$-valued functions on $\Omega$ and the point evaluations in 
$t\in\Omega$ are continuous linear functionals on both $\BB$ and $\BB^\#$.
\end{deff}

Since we want to use the results of the previous section, we are going to work with a special class of \RKBSS. We call a \UU\ \RKBS\ by a \textit{\SIP\ \RKBS} (\SIPRKBS\ for short). While it is true that a RKBS possesses some sort of function that resembles to the \RK\ 
for a RKHS, in a \SIPRKBS\ we have a function with those same attributes of the \RK\ for a RKHS. 

\begin{prop}\label{rk} 
Let $\BB$ be a RKBS on $\Omega$, then there exists a unique function (\RK) $K:\Omega\times\Omega\rightarrow\CC$ such that:

(1) For all $t\in\Omega$, $K(\cdot , t)\in\BB^*$ and $f(t) = (f,K(\cdot , t))_{\BB}$ for all $f\in\BB$.

(2) For all $t\in\Omega$, $K(t,\cdot )\in\BB$ and $f^*(t) = (K(t, \cdot ),f^*)_{\BB}$ for all $f^*\in\BB^*$.  

(3) $\BB^* = \overline{span}\{ K(\cdot , t) : t\in\Omega \}$ and $\BB   = \overline{span}\{ K(t, \cdot ) : t\in\Omega \}$.

(4) $K(s,t)= (K(s, \cdot ),K(\cdot , t))_{\BB}$ for all $s,t\in\Omega$.
\\
Moreover, if $\BB$ is also a \SIPRKBS\ on $\Omega$, then there exists another unique function (\SIPK) 
$G:\Omega\times\Omega\rightarrow \CC$ such that

(5) $ G(t, \cdot ) \in\BB$ and $K(\cdot ,t)= (G(t, \cdot ))^* \in\BB^*$ for all $t\in\Omega$.

(6) $f(t)= [f,G(t,\cdot )]$ and $f^* (t)= [K(t,\cdot ),f]$ for all $f\in\BB ,\, t\in\Omega$.
\\
When $K=G$, we call it the \SIPRK\ for $\BB$.
\end{prop}

An important result in a RKHS is that norm convergence implies pointwise convergence, the same is true in a RKBS (therefore in a \SIPRKBS).
An another one is about how it can be constructed a \SIPRKBS\ by using an isometric operator. The following construction appears in 
\cite{Gar2, Zhang1}. 

\begin{comm}[ \SIPRKBS\ construction by using an operator]\label{operator-construction}  
Let $(E , [ \cdot ,\cdot ]_E)$ be a \UU\ \BS; let $\Phi :\Omega\rightarrow E$ be a function and let $T_\Phi : E \rightarrow \CC^\Omega$ be 
an operator defined by $T_\Phi x=f_x$ with $f_x (t) = [ x, \Phi (t)]_E$, $t\in\Omega$. 
It follows that $T_\Phi$ is linear, and it is injective if we suppose further that $ \{ \Phi (t):t\in\Omega \}$ is a complete set in $E$. 
Let $\BB= \rr(T_\Phi )$ be the range of $T_\Phi$ and define the $\BB$-norm by $\Vert f_x \Vert_{\BB} := \Vert x \Vert_E$, this turns $T_\Phi$ into an isometric 
isomorphism between $E$ and $\BB$, therefore $\BB$ is a \UU\ \BS\ of $\CC$-valued functions on $\Omega$. 
Moreover, $[ \cdot ,\cdot ]_{\BB}$ defined by $[f_x , f_y]_{\BB} :=[x,y]_E \ x,y\in E$ is the unique (norm compatible) s.i.p. on $\BB$. 
For each $t\in\Omega$ the point evaluations over $\BB$ are continuous but, for being continuous over $\BB^*$ we need some extra hypotheses. 
We consider the function $\Phi^* :\Omega \rightarrow E^*$ given by $\Phi^* (t) =(\Phi (t))^*,\ t\in\Omega$, and impose that 
$\overline{span} \{ \Phi^* (t):t\in\Omega \} = E^*$. 
In this way (see \cite[Theo. 10]{Zhang1}) $\BB^* = \{ f_x^* :=[ \Phi (\cdot ), x]_{E} : x\in E \}$ endowed with 
$[ f_x^* , f_y^* ]_{\BB^*} := [f_y , f_x ]_{\BB}\ \, x,y \in E$ is the dual of $\BB$ with the bilinear form $( f_x , f_y^* )_{\BB} := (x,y^*)_E \ \, x,y\in E$ which 
\SIPRK\ $G$ for $\BB$ is given by $G(s,t) = [\Phi (s) , \Phi (t)]_{E} \ \, s,t\in\Omega$.
\end{comm}

If it is necessary to distinguish each characteristic component of a \SIPRKBS\ on $\Omega$ constructed as before, then we write it as 
$(\BB , [ \cdot ,\cdot ]_{\BB}, G, E, \Phi)$. 

In the spirit of Zayed's book \cite[Def. 10.1.3.]{Zayed} we have the following definition which is our main objective when we talk about reconstruction 
in sampling theory.

\begin{deff}[Sampling Basis]
A basis $\{ S_j \}_{j\in\II}$ of a \RKBS\ $\BB$ on a subset $\Omega$ is called a \textit{sampling basis} if there exists a sequence 
$\{t_j\}_{j\in\II} \subset\Omega$ such that
\begin{equation}\label{eq:samplingbasis}
                f(t)=\sum_{j\in\II} f(t_j)S_j (t)
                \qquad\forall\,f\in\BB ,\ t\in\Omega
\end{equation} 
\end{deff}

Again, by similarity with the \HS\ sampling theory, we need to restrict to work in a \SIPRKBS.

\begin{prop}[Sampling basis]\label{samplingbasisiff}
Let $\{ S_j \}_{j\in\II}$ be a basis of a \SIPRKBS\ $\BB$ on a set $\Omega$ with \SIPRK\ $G$. Then, $\{ S_j \}_{j\in\II}$ is a 
sampling basis \IFF\ its \BIO\ basis $\{ F_j \}_{j\in\II}$ is given by
\begin{equation}
                F_j(t)= G(t_j , t):=G_{t_j}(t)
                \qquad j\in\II ,\ t\in\Omega
\end{equation}
\begin{proof}  
$( \Rightarrow )$ Since $\{S_{j}\}_{j\in\II}$ is a sampling basis for $\BB$ its (unique) \BIO\ \SB, for each $f\in\BB$, 
$\{ F_j \}_{j\in\II}$ satisfies:
\begin{displaymath}
                   \sum_{j\in\II} [ f, F_j ] S_{j}(t)
                   = f(t)   
                   = \sum_{j\in\II} f(t_j) S_{j}(t)  
                   = \sum_{j\in\II} [ f, G_{t_j} ] S_{j}(t)  
                   \qquad t\in\Omega
\end{displaymath}
In the first equality is used the \SB\ property of both sequences, while the second equality follows by the sampling basis hypothesis over 
$\{S_{j}\}_{j\in\II}$ and the last one is due to the reproducing property of $G$. Thus, by uniqueness, it must be 
$[ f, F_j ] = [ f, G_{t_j} ] \ \,\forall\,j\in\II$, whence $F_j = G_{t_j}\ \,\forall\,j\in\II$.
             
$( \Leftarrow )$ If the \BIO\  \textnormal{\SB} to $\{S_{j}\}_{j\in\II}$ is given by $F_j (t)= G(t_j,t):= G_{t_j} (t)$, then for each 
$f\in\BB$:
\begin{displaymath}
                   f(t) = \sum_{j\in\II} [ f, F_j ] S_{j}(t) = \sum_{j\in\II} [ f, G_{t_j} ] S_{j}(t)  
                        = \sum_{j\in\II} f(t_j) S_{j}(t)   
                   \qquad t\in\Omega
\end{displaymath}
therefore $\{S_{j}\}_{j\in\II}$ is a sampling basis for $\BB$.
\end{proof} 
\end{prop}
Of course, the definition of sampling basis as well as the last proposition are valid in a RKHS because every RKHS is a \SIPRKBS.
%
\section{Kramer-Type Sampling Theorems}\label{sec:krasampling}
\markboth{\thesection. Kramer-Type Sampling Theorems}{\thesection. Kramer-Type Sampling Theorems}
%
The next procedure for obtaining a \SIPRKBS\ version of the \KST\ is due to Garc\'ia, Hern\'andez-Medina \& Mu\~{n}oz-Bouzo \cite{Gar3}, 
they use a BK-space instead of the $\ell_2$ space, an $X_d^*$-\textit{\RB} instead of a \RB\ and a s.i.p. RKBS instead of a RKHS. 
Let $(E,[ \cdot ,\cdot ]_E)$, $\Phi :\Omega\rightarrow E$ and $T_\Phi :E\rightarrow\CC^\Omega$ be as in Remark \ref{operator-construction}. 
First we suppose there exists a sequence $\{x_{j}\}_{j\in\II}\subset E$ such that $\{x_{j}^*\}_{j\in\II}$ is an $X_d^*$-\textit{\RB} for 
$E^*$, then, there exists an unique \BIO\ sequence $\{ y_j \}_{j\in\II}$ which is an $X_d$-\textit{\RB} for $E$ (see \cite{Zhang1}). In 
second place, suppose the existence of sequences $\{ t_j \}_{j\in\II} \subset\Omega$ and $\{ a_j \}_{j\in\II} \subset \CC \setminus \{ 0 \}$ 
such that the interpolation condition $(\Phi (t_k) )^* = a_k x_k^*\ k\in\II$ or, equivalently, 
$S_j (t_k ): =[y_j , \Phi (t_k)]_E= a_j \delta_{j,k}\ j,k\in\II$ holds true, where for fixed $t\in\Omega$ we have
$( \Phi (t) )^* = \sum_{j\in\II} S_j (t)\, x_j^* \,\in E^*$ with $S_j (t):=[y_j , \Phi (t)]_E$, being the sequence 
$\{S_j(t)\}_{j\in\II}\in X_d^*$ for each $t\in\Omega$ as can be checked by strightforward calculations. 
Under these hypotheses is obtained the \SIPRKBS\ $\BB$ on $\Omega$ explicitly given by 
$\BB =\{ f_x(\cdot )=[ x, \Phi(\cdot )]_E: x\in E\}$, with norm $\Vert f_x\Vert_{\BB} :=\Vert x\Vert_E$ and \SIPRK\ 
$G(s,t) = [\Phi (s) , \Phi (t)]_E \ s,t\in\Omega $. We now state the before mentioned {\SIPRKBS} version of the \KST.

\begin{theo}[\SIPRKBS\ \KST\ {\cite[p. 19]{Gar3}}]\label{kstsiprkbs}
Let $\BB$ be a \SIPRKBS\ on $\Omega$ as before. Then, the sequence  $\{S_j\}_{j\in\II} \subset \BB$ is an 
$X_d$-\textit{\RB} for $\BB$ and for any $f\in\BB$ we have the sampling expansion
\begin{displaymath}
                   f(t) = \sum_{j\in\II} f(t_j) \frac{S_j (t)}{a_j} 
                   \qquad t\in\Omega
\end{displaymath}
The series converges in the $\BB$-norm sense and also, absolutely and uniformly on subsets of $\Omega$ where the function 
$t \mapsto \Vert \Phi (t) \Vert_E$ is bounded.
\end{theo}

\begin{coro}
Under hypotheses of Theorem \ref{kstsiprkbs}, $\{a_j^{-1}S_j\}_{j\in\II}$ is a sampling basis for $\BB$.
\begin{proof}  
Indeed, for each $j,k\in\II$ we have
\\

$\displaystyle \qquad [a_j^{-1}S_j , G_{t_k}]_{\BB} = \bigg[\frac{S_j}{a_j}  , G_{t_k} \bigg]_{\BB} = 
               \frac{1}{a_j}[S_j , G_{t_k}]_{\BB} = \frac{1}{a_j}S_j (t_k) = \frac{a_j}{a_j}\delta_{j,k} = \delta_{j,k}$.
\end{proof}
\end{coro}
%
\section{The main result: A converse of the Kramer Sampling Theorem in a \SIPRKBS}\label{sec:mainresult}
\markboth{\thesection. The main result: A converse of the \KST\ in \SIPRKBS}
{\thesection. The main result: A converse of the \KST\ in \SIPRKBS}
%
Keeping in mind the statement as well as the proof of the Theorem \ref{conversekst} (\cite[pp. 55 --58]{Gar1}), we consider 
$(\BB , [ \cdot ,\cdot ]_{\BB}, G, E, \Phi)$ a \textnormal{\SIPRKBS} on $\Omega$ that it has been built by an isometric operator (Remark 
\ref{operator-construction}) and we assume the existence of a sequence $\{ S_j \}_{j\in\II}\subset \BB $ such that
\begin{equation}\label{eq:Sj-xd}
                \{ S_j (t) \}_{j\in\II} \in X_d 
                \qquad \textrm{and} \qquad
                \{ S_j (t) \}_{j\in\II}^* \in X_d^* \qquad  \forall\, t\in\Omega   \nonumber    
\end{equation}
with $X_d$ a \UU\ \textit{BK-space} (for instance $\ell_p (\II)$). We define two functions 
\begin{displaymath}
                   \begin{array}{ccc}
                   \phi :\Omega & \longrightarrow & X_d \qquad\ \\
                   \quad\ t     & \longmapsto     & \{ S_j(t) \}_{j\in\II}
                   \end{array}
                              \qquad \textnormal{ and } \qquad
                   \begin{array}{ccc}
                   \phi^* :\Omega & \longrightarrow & X_d^* \qquad\ \ \\
                   \qquad t       & \longmapsto     & \{ S_j(t) \}^*_{j\in\II}
                   \end{array}
\end{displaymath}
using the notation $\phi^* (t)= (\phi (t))^*$, $t\in\Omega$. We also assume that holds true:
\begin{equation}\label{eq:Sj-completas}
\left\lbrace \begin{array}{ccccccc}
\textnormal{if } & c\in X_d & \textnormal{ is such that} & \sum_{j\in\II} c_j(S_j (t))^*=0 & \forall\,t\in\Omega & \Rightarrow & c=0. \\
\textnormal{if } & d\in X_d^* & \textnormal{ is such that} & \sum_{j\in\II} d_j S_j (t)=0 & \forall\, t\in\Omega & \Rightarrow & d=0.
\end{array} \right.
\end{equation}
This requirement is similar to that in the item $1^\circ )$ of Theorem \ref{conversekst}, and it is equivalent to the completeness statement:
\begin{equation}
                   \overline{span} \{ \phi(t):t\in\Omega \}=X_d 
                   \qquad \textrm{and} \qquad
                   \overline{span} \{ \phi^*(t):t\in\Omega \}=X_d^* 
\end{equation}
which is necessary for the definition itself of the \textnormal{\SIPRKBS} 
$(\BB_{samp} , [ \cdot ,\cdot ]_{samp}, G_{samp}, X_d^*, \phi^*)$ on $\Omega$. By the way, its \textnormal{\SIPRK} $G_{samp}$ is given by 
\begin{equation}\label{eq:G-samp}
                G_{samp}(s,t):= [ \phi^* (s), \phi^* (t) ]_{X_d^*} 
                = \sum_{j\in\II} (S_j (s))^* S_j (t)
                \qquad s,t\in\Omega 
\end{equation}
where the reflexivity of $X_d$ was used to the identification of $(\phi (t))^{**}$ with $\phi(t)$. 

We have taken $X_d^*$ instead of $X_d$ in the definition of $\BB_{samp}$ because we want the similarity between the 
\textnormal{\SIPRK} $G_{samp}$ and the \textnormal{\RK} $K_{samp}$, where the last one was used in Theorem \ref{conversekst}.
We are going to prove three propositions that will be used in the demonstration of the main result, the first two are interesting on 
their own.

\begin{prop}\label{Sj-bessel}
Let $(\BB , [ \cdot ,\cdot ]_{\BB}, G, E, \Phi)$ and $(\BB_{samp} , [ \cdot ,\cdot ]_{samp}, G_{samp}, X_d^*, \phi^*)$ 
be two \SIPRKBS\ on $\Omega$ as before. If the sets $\{\, \{S_j (t)\}_{j\in\II} : t\in\Omega \}\subset  X_d$ and 
$\{\, \{S_j (t)\}^*_{j\in\II} : t\in\Omega \}\subset X_d^*$ are complete, then $\{S_j^* \}_{j\in\II}$ is an $X_d$-\textit{Bessel} 
sequence for $\BB^*$. 
\begin{proof}  
The completeness conditions (equivalent to \eqref{eq:Sj-completas}) are stated because it is necessary for the definition of $\BB_{samp}$. 
We must to show there exist $B>0$ such that 
\begin{displaymath}
                  \bigg\Vert \sum_{j\in\II} d_j S_j \bigg\Vert_{\BB} \leq B\Vert d\Vert_{X_d^*} 
                  \qquad\forall\, d\in X_d^* 
\end{displaymath}
For Proposition \ref{besselseqiff} it is equivalent to show the associated analysis operator given by 
\begin{displaymath}
                   \begin{array}{ccc}
                   V:\BB^*&\longrightarrow & X_d\quad\quad\quad\quad\;\; \\
                   \quad\; f^* &\longmapsto & \{[ f^*,S_j^* ]_{\BB^*} \}_{j\in\II}
                   \end{array}
\end{displaymath}
is bounded. The operator $T: E^*\rightarrow\BB^*$ defined by $Tx^*=[x^*,\Phi^* (\cdot )]_{E^*}:=f_{x^*}$ is an isometric isomorphism, 
therefore it sends dense subspaces on $E^*$ in dense subspaces on $\BB^*$. We know it suffices to prove the {Bessel} condition of 
$\{S_j^* \}_{j\in\II}$ on a dense subset of $\BB^*$. Since the set $span \{ \Phi^* (s) : s\in \Omega \}$ is dense in $E^*$ and
\begin{equation}
                T \Phi^* (s) = [\Phi^* (s), \Phi^* (\cdot ) ]_{E^*} = [\Phi (\cdot )  , \Phi(s) ) ]_{E} = 
                G(\cdot , s) = (G(s, \cdot ))^* := G_s^* 
                \qquad s\in\Omega \nonumber
\end{equation}
the set $\BB^*_0=span \{ G_s^* : s\in \Omega \}$ is dense in $\BB^*$. Now, we consider for each $N\in\NN$
\begin{displaymath}
                   \begin{array}{ccc}
                   V_N :\BB^*_0 &\longrightarrow & X_d \qquad\qquad\qquad\quad\; \\
                   \qquad f^* &\longmapsto &\{ \textbf{1}_{\II_{N}}(j) [ f^*,S_j^* ]_{\BB^*} \}_{j\in\II}
                   \end{array}
                              \quad \textnormal{ and } \quad
                   \begin{array}{ccc}
                   V' :\BB^*_0 &\longrightarrow & X_d \qquad\qquad\;\; \\
                   \quad\;\; f^* &\longmapsto &\{[ f^*,S_j^* ]_{\BB^*} \}_{j\in\II}
                   \end{array}
\end{displaymath}
where $\textbf{1}_{\II_{N}}$ denotes the characteristic function of $\II_{N}$ (the first $N$ elements of $\II$). For each $s\in\Omega$, 
$j\in\II$ there holds:
\begin{displaymath}
                   [G^*_s , S_j^* ]_{\BB^*} = [ S_j , G_s]_{\BB} = S_j(s)
\end{displaymath}
and since $\{ S_j (s)\}_{j\in\II} \in X_d\ \,\forall\, s\in\Omega$, the operators $V_N$ are well defined and also they are bounded for 
each $N\in\NN$, since
\begin{displaymath}
                   \Vert V_N f^* \Vert_{X_d} =
                   \sup_{d\in S_{X_d^*}} \bigg \vert \sum_{j\in\II_N} d_j [ f^*,S_j^* ]_{\BB^*} \bigg \vert \leq
                   \sup_{d\in S_{X_d^*}}\bigg ( \sum_{j\in\II_N}  d_j  \Vert  S^*_j \Vert_{\BB^*} \bigg ) \Vert  f^* \Vert_{\BB^*}
\end{displaymath}
Furthermore, they converge pointwise to V' since
\begin{align*}
              \Vert V_N f^* - V'f^* \Vert_{X_d} =
              \sup_{d\in S_{X_d^*}} \bigg \vert \sum_{j\in \II \setminus \II_N} d_j [ f^*,S_j^* ]_{\BB^*} \bigg \vert \leq
              \Vert \{ \textbf{1}_{\II \setminus \II_{N}}(j)[ f^*,  S_j^* ]_{\BB^*} \}_{j\in \II} \Vert_{X_d}
              \; \mathop{\quad\longrightarrow 0}_{N\rightarrow\infty\;} 
\end{align*}
Thus, by \textit{Banach-Steinhaus theorem}, $V'$ is a bounded operator, therefore $V$ it is, and $\{S_j^* \}_{j\in\II}$ is an $X_d$-\textit{Bessel} sequence for $\BB^*$. 
\end{proof}
\end{prop}

An infinite-dimensional \VS\ can be endowed with various norms which turns it in a \BS, but being non-equivalent between them (by the existence of unbounded linear 
functionals). Of course, this phenomenon does not occur in a finite-dimensional \BS\ and we prove in the following proposition that neither occurs in a \RKBS, due to the convergence property: if $f_j$ converges to $f$ in $\BB$, then $f_j$ converges pointwise to $f$ in $\Omega$ \cite{Zhang1}.

\begin{prop}\label{normequivalent}
Let's suppose that $\BB$ is a \SIPRKBS\ on $\Omega$ endowed with the norm $\Vert \cdot\Vert_{\BB}$ either the norm $\Vert \cdot\Vert$. Then, the norms are equivalent.
\begin{proof}  
We show the identity operator $id:(\BB , \Vert \cdot\Vert_{\BB})\rightarrow (\BB , \Vert \cdot\Vert)$ is bounded, and then by the 
\textsl{open mapping theorem} will result bi-continuous. By the \textsl{closed graph theorem}, we only need to check:
\begin{displaymath}
                   f_j \longrightarrow f 
                         \quad \textnormal{in}\quad
                   \Vert \cdot \Vert_{\BB}
                         \qquad \textnormal{and}\qquad
                   f_j \longrightarrow g
                         \quad \textnormal{in}\quad
                   \Vert \cdot \Vert
                         \qquad \textnormal{then}\qquad
                   f=g
\end{displaymath}
and this is clear due to the convergence property in a \SIPRKBS.
\end{proof}
\end{prop}

\begin{prop}\label{Sj-rieszbasis}
Let $(\BB , [ \cdot ,\cdot ]_{\BB}, G, E, \Phi)$ and $(\BB_{samp} , [ \cdot ,\cdot ]_{samp}, G_{samp}, X_d^*, \phi^*)$ be two \SIPRKBS\ 
on $\Omega$ as before. Let's suppose that:

\begin{itemize}
\item[$1^\circ$)] The sets $\{\, \{S_j (t)\}_{j\in\II} : t\in\Omega \}\subset  X_d$ and 
                  $\{\, \{S_j (t)\}^*_{j\in\II} : t\in\Omega \}\subset X_d^*$ are complete.

\item[$2^\circ$)] There exists sequences $\{t_j\}_{j\in\II} \subset\Omega ,\; \{a_j\}_{j\in\II} \subset \CC \setminus \{ 0 \}$ such 
                  that there holds the following sampling conditions:
                  
                  \begin{equation}\label{eq:formula-sampling-1}
                  \bigg \{ \frac{f(t_j )}{a_j} \bigg \}_{j\in\II} \in X_d^*
                  \qquad\forall\, f\in\BB
                  \end{equation}
                  and
                  \begin{equation}\label{eq:formula-sampling-2}
                  f(t) = \sum_{j\in\II} f(t_j)\frac{S_j(t)}{a_j}
                  \qquad\forall\, f\in\BB
                  \end{equation}
                  where the series converges absolutely on $\Omega$. 
\end{itemize}
If we call: 
\begin{equation}\label{eq:notacion-Mj-gsamp}
                M_j(\cdot):=\overline{a_j}^{-1}G_{samp}(t_j,\cdot )
                           \qquad \textnormal { and } \qquad
                M_j^* = a_j^{-1}G_{samp}(\cdot ,t_j )\in\BB_{samp}^*
                           \qquad j\in\II 
\end{equation}
Then:
\begin{itemize}
               \item[a)] $\{ M_j^*\}_{j\in\II}$ is a complete sequence in $\BB_{samp}^*$.
               \item[b)] $\{ M_j \}_{j\in\II}$ is an $X_d^*$-\textit{Bessel} sequence for $\BB_{samp}$.
               \item[c)] $\{ M_j\}_{j\in\II}$ is a minimal sequence in $\BB_{samp}^*$ with \BIO\ sequence $\{S_j\}_{j\in\II}$. Also, 
                         $\{ M_j^*\}_{j\in\II}$ is a minimal sequence in $\BB_{samp}$ with \BIO\ sequence $\{S_j^*\}_{j\in\II}$.
               \item[d)] $\{ M_j\}_{j\in\II}$ is an $X_d^*$-\textit{Riesz-Fischer} sequence for $\BB_{samp}$.
               \item[e)] $\{ M_j\}_{j\in\II}$ is an $X_d^*$-\textit{frame} for $\BB_{samp}$.
               \item[f)] $\{ M_j\}_{j\in\II}$ is an $X_d^*$-\textit{\RB} for $\BB_{samp}$.        
\end{itemize}
\begin{proof}  
$a$) We assume there exists $f\in\BB_{samp}$ such that $0 = [f, M_j]_{samp}=[M_j^*,f^*]_{samp*}\ \forall\, j\in\II$. But, since
\begin{displaymath}
                   f(t) = \sum_{j\in\II} f(t_j) \frac{S_j (t)}{a_j} = \sum_{j\in\II} [f, M_j]_{samp} \,S_j (t) = 0 
                   \qquad \forall\, t\in\Omega 
\end{displaymath}
then $f=0$.
             
${b}$) This is immediate since $\{ f(t_j) a_j^{-1}  \}_{j\in\II} \in X_d^*\ \,\forall \,f\in\BB_{samp}$ and, $\forall \, j\in\II$, holds
\begin{align*} 
              a_j^{-1}f(t_j) =
              a_j^{-1} [f, G_{samp} (t_j , \cdot )]_{samp} =
              [f,\overline{a_j}^{-1} G_{samp} (t_j , \cdot )]_{samp} =
              [f,M_j]_{samp}
\end{align*}
Because of this we also obtain the well-definition and boundedness of the analysis operator $U:\BB_{samp} \rightarrow X_d^*$ 
associated to the sequence $\{ M_j\}_{j\in\II}$ as well as its adjoint $U^* :X_d \rightarrow \BB_{samp}^*$, in particular 
$U^* c= \sum_{j\in\II} c_j M_j^*$ converges (unconditionally) in $\BB_{samp}^*$ for all $c\in X_d$.

$c$) Due to $[ M^*_k , S^*_j ]_{samp*}=[ S_j , M_k ]_{samp}$, we only need to show that 
$[ S_j , M_k ]_{samp} = \delta_{j,k}\ \forall\, j,k\in\II$. In one hand we have
\begin{displaymath}
                   S_k(t) = \sum_{j\in\II} \delta_{j,k} S_j(t) 
                   \qquad k\in\II ,\,t\in\Omega
\end{displaymath}
and by other hand
\begin{displaymath}
                   S_k(t) = \sum_{j\in\II} \frac{S_k (t_j)}{a_j} S_j(t) = \sum_{j\in\II} [ S_j , M_k ]_{samp} \,S_j(t) 
                   \qquad k\in\II ,\,t\in\Omega
\end{displaymath}
then
\begin{displaymath}
                   0 = \sum_{j\in\II} \big ( [ S_j , M_k ]_{samp} - \delta_{j,k} \big )S_j(t)
                   \qquad k\in\II ,\,t\in\Omega
\end{displaymath}
where the coefficients are in $X_d^*$, therefore we obtain $[ S_j , M_k ]_{samp} = \delta_{j,k}\ \,\forall\, j,k\in\II$ as we needed.

$d$) Given a sequence $d=\{d_j \}_{j\in\II} \in X_d^*$ we must to see there exist $f\in\BB_{samp}$ such that $Uf=d$. By considering 
$f=\sum_{j\in\II} d_j S_j$ (it belongs to $\BB_{samp}$) it leads to
\begin{align*}
              \big [ f, M_k \big ]_{samp} = \bigg [ \sum_{j\in\II} d_j S_j, M_k \bigg ]_{samp} 
              = \sum_{j\in\II} d_j \big [ S_j, M_k \big ]_{samp} = d_k
              \qquad k\in\II
\end{align*}
therefore $\{Uf\}_k =d_k \ \,\forall\, k\in\II$ and $U$ is surjective.

$e$) It follows by items $b$) and $d$) due to Proposition \ref{framesiff}.

$f$) It follows by items $c$) and $e$) due to Proposition \ref{rieszbasisiff}.
\end{proof}           
\end{prop}

We now prove the main result of this paper.

\begin{theo}[A Converse of the \KST\ - \SIPRKBS\ Version]\label{recipkstsiprkbs}
Under hypotheses of Proposition \ref{Sj-rieszbasis} we have:
                  \begin{itemize}
                  \item[a)] $\BB_{samp} = \BB$.
                                                                    
                  \item[b)] The norms $\Vert \cdot \Vert_{\BB_{samp}}$ and $\Vert \cdot \Vert_{\BB}$ are equivalent and 
                            consequently $\{S_j\}_{j\in\II}$ is an $X_d^*$-\textit{\RB} for $\BB$.
                                                                          
                  \item[c)] The \BIO\ sequence of $\{ S_j\}_{j\in\II}$ in $\BB_{samp}$ is given by
                  \begin{equation}\label{eq:formula-sampling-3}
                                  \bigg \{
                                               \sum_{k\in\II}\bigg[ \frac{\phi^*(t_j)}{\overline{a_j}},
                                               \frac{\phi^*(t_k)}{\overline{a_k}}\bigg ]_{X_d^*} \!\! S_k
                                  \bigg \}_{j\in\II}   
                  \end{equation}
                                                              
                  \item[d)] The \BIO\ sequence of $\{ S_j\}_{j\in\II}$ in $\BB$ is given by
                  \begin{equation}\label{eq:formula-sampling-4}
                                  \bigg \{ 
                                               \sum_{k\in\II}\bigg [\frac{\Phi(t_j)}{\overline{a_j}},
                                               \frac{\Phi(t_k)}{\overline{a_k}}\bigg ]_{E} \!S_k
                                  \bigg \}_{j\in\II}   
                  \end{equation}
                  \end{itemize}                                                        
\begin{proof}
$a$) We first prove that $\BB_{samp} \subset \BB$ by only assuming the item $1^\circ)$. Due to $\BB_{samp}$ comprises functions of 
the form
\begin{equation}
                \sum_{j\in\II} d_j S_j
                \quad \textnormal{ with } \quad 
                d=\{ d_j \}_{j\in\II} \in X_d^* \nonumber
\end{equation}
by definition, it follows that $\Vert \sum_{j\in\II} d_j S_j \Vert_{samp} < \infty \ \,\forall\, d\in X_d^*$. To see 
$\BB_{samp}\subset\BB$ we must to show $\Vert \sum_{j\in\II} d_j S_j \Vert_{\BB} < \infty \ \,\forall\, d\in X_d^*$. By Proposition 
\ref{Sj-bessel} $\{S_j^* \}_{j\in\II}$ is an $X_d$-\textit{Bessel} sequence for $\BB^*$, therefore the analysis operator associated 
to $\{S_j^* \}_{j\in\II}$ is bounded and so it is the synthesis operator, i.e.:
\begin{displaymath}
                   \bigg\Vert \sum_{j\in\II} d_j S_j \bigg\Vert_{\BB} \leq B \Vert d \Vert_{X_d} < \infty
                   \qquad \textnormal{for some $B>0$}
\end{displaymath}
             
For the other inclusion we also assume to hold true the sampling conditions \eqref{eq:formula-sampling-1} and 
\eqref{eq:formula-sampling-2}. We pick $f\in\BB$, then $\{ f(t_j) a_j^{-1}  \}_{j\in\II} \in X_d^*$ by \eqref{eq:formula-sampling-1} and 
the series $\sum_{j\in\II} f(t_j) a_j^{-1}S_j$ converges in $\Vert \cdot \Vert_{samp}$, we say to $g\in\BB_{samp}$, therefore it 
converges pointwise to $g\in\BB_{samp}$, but the series also converges pointwise to $f$ by \eqref{eq:formula-sampling-2}, whence 
$g(t)=f(t)\ \,\forall\, t\in\Omega$ and hence $f\in\BB_{samp}$.
             
$b$) As we have $\BB =\BB_{samp}$, the equivalence between the norms $\Vert \cdot \Vert_{\BB}$ and $\Vert \cdot \Vert_{samp}$ follows by 
Proposition \ref{normequivalent} and since $\{S_j\}_{j\in\II}$ is the \BIO\ sequence to $\{ M_j\}_{j\in\II}$ 
(Proposition \ref{Sj-rieszbasis}, item c)), it is an $X_d$-\textit{\RB} for $\BB^*_{samp}$ \cite[Theo. 2.14 and 2.15]{Zhang2} 
as well as an $X_d$-\textit{\RB} for $\BB^*$ by norm equivalence. 
             
We recall the notations \eqref{eq:notacion-Mj-gsamp} and now we add a new one: $ G_j(\cdot):=\overline{a_j}^{-1}G(t_j,\cdot )\ \,j\in\II$.
             
$c$) We have already seen that $\{ S_j \}_{j\in\II}$ and $\{ M_j\}_{j\in\II}$ are \BIO\ sequences in $\BB_{samp}$, so we are going to see 
there holds \eqref{eq:formula-sampling-3}, indeed for $k\in\II ,\,t\in\Omega$ we have
\begin{align*}
              M_k(t) =
              \sum_{j\in\II}\frac{M_k(t_j)}{{a_j}} S_j (t) =
              \sum_{j\in\II}[M_k,M_j]_{samp} \,S_j(t) =
              \sum_{j\in\II}\bigg[ \frac{\phi^*(t_k)}{\overline{a_k}},\frac{\phi^*(t_j)}{\overline{a_j}}\bigg ]_{X_d^*} \!\!S_j (t)
\end{align*}
                       
$d$) Again, we have already seen that $S_j(t_k)=a_k \delta_{j,k}\ \,\forall\,j,k\in\II$, whence 
\begin{displaymath}
                   \delta_{j,k} = \frac{S_j (t_k)}{a_k} = \bigg [ S_j , \frac{G_{t_k}}{\overline{a_k}} \bigg ]_{\BB} = [ S_j , G_k ]_{\BB} 
                   \qquad j,k\in\II                                       
\end{displaymath}
and therefore $\{ S_j \}_{j\in\II}$ and $\{ G_j \}_{j\in\II}$ are \BIO\ sequences in $\BB$. Finally, $\{ G_j \}_{j\in\II}$ satisfies 
\eqref{eq:formula-sampling-4}, since
\begin{align*}
              G_k (t) =
              \sum_{j\in\II}\frac{G_k(t_j)}{{a_j}} S_j (t) =
              \sum_{j\in\II} [ G_k  , G_j ]_{\BB} \,S_j (t) =
              \sum_{j\in\II}\bigg [\frac{\Phi(t_k)}{\overline{a_k}},\frac{\Phi(t_j)}{\overline{a_j}}\bigg ]_{E} \!S_j(t)
\end{align*}
for all $k\in\II ,\,t\in\Omega$. This finishes the proof.
\end{proof}
\end{theo}

\begin{coro}
Under hypotheses of Theorem \ref{recipkstsiprkbs}, $\{a_j^{-1}S_j\}_{j\in\II}$ is a sampling basis for $\BB_{samp}$.
\begin{proof}  
By Proposition \ref{samplingbasisiff} we only need to check $[a_j^{-1}S_j , G_{t_k}]_{samp}=\delta_{j,k}\ \forall\,j,k\in\II$ since 
$\{ a_j^{-1}S_j\}_{j\in\II}$ is a \SB. We have
\begin{equation}
                [a_j^{-1}S_j , G_{t_k}]_{samp}=
                \bigg[\frac{S_j}{a_j}  , \overline{a_k}\frac{G_{t_k}}{\overline{a_k}} \bigg]_{samp}=
                \frac{a_k}{a_j}[S_j , M_k]_{samp}=
                \frac{a_k}{a_j}\delta_{j,k}=
                \delta_{j,k} \nonumber
\end{equation}
as we needed.
\end{proof}
\end{coro}

We finish with a classical example.

\begin{exam}
We consider $\Omega=\RR$, $I=\left[-\frac{1}{2},\frac{1}{2}\right]$, $1\!<\!p, q\!<\! +\infty$ with $\frac{1}{p}+\frac{1}{q}=1$, the ``{Time-limited}'' \UU\ \BSS: 
\begin{eqnarray}
                L_p(I,\RR)&:=&\{ f\in L_p(\RR): f\equiv 0 \textnormal{ a.e. on }\RR\setminus I\} \nonumber
\end{eqnarray}
and the ``{band-limited}'' \UU\ \BSS:
\begin{eqnarray}
                \BB_p:&=&\{ f\in C(\RR )\cap L_q (\RR): \widehat{f}\equiv 0 \textnormal{ a.e. on }\RR\setminus I  \} \nonumber
\end{eqnarray}

We define $\Phi :\Omega\rightarrow L_p(\RR)$ and $\Phi^* :\Omega\rightarrow L_q(\RR)$ by 
\begin{displaymath}
                  \Phi (\omega)(t ):= e^{-2\pi it\omega }
                               \qquad \textnormal{and} \qquad
                  \Phi^* (\omega)(t ):= e^{2\pi it\omega }
                               \qquad \omega\in \RR,\; t\in \RR
\end{displaymath}           
It is well-known that
\begin{eqnarray}
                \overline{span} \{ \chi_I(\cdot ) e^{-2\pi i(\cdot )\omega } : \omega\in \RR \} &=& L_p(I) \nonumber
\end{eqnarray}
            
Let $\TF$ be the Fourier transform, $\TF:L_1(\RR)\rightarrow C_0 (\RR)$ given by
\begin{equation}
                \widehat{f}(\omega ):= \TF [f](\omega )= \int_{\RR}f(t)e^{-2\pi it\omega}dx \nonumber
\end{equation}
and we note by $f^{\vee}$ the Fourier inversion of $f$ given by
\begin{equation}
                f^{\vee}(\omega ):= \TF^{-1} [f](\omega )= \int_{\RR}f(t)e^{2\pi it\omega}dt\nonumber
\end{equation}
Clearly, $\widehat{f}(\omega )=f^{\vee} (-\omega )\ \,\forall\,\omega\in\RR$, and by Fourier Analysis we know that if $f, \widehat{f}\in L_1(\RR )$ then $f$ and 
$\widehat{f}$ are continuous and we have the inversion formulae:
\begin{eqnarray}
                f(t)= &\widehat{({f}^{\vee})}(t) =&\int_{\RR}f^{\vee}(\omega)e^{-2\pi it\omega}d\omega     \nonumber \\
                f(t)= &(\widehat{f}\,)^{\vee}(t) =&\int_{\RR}\widehat{f}\,(\omega)e^{2\pi it\omega}d\omega \nonumber
\end{eqnarray} 
where the equality is pointwise $t\in\RR$. Also there holds
\begin{eqnarray}
                \int_{\RR}f(\omega)\widehat{g}(\omega)d\omega        &=& \int_{\RR}\widehat{f}(\omega) g(\omega)d\omega \nonumber \\
                \int_{\RR}\widehat{f}(\omega)g^{\vee}(\omega)d\omega &=& \int_{\RR}f(x){g}(x)dx  .        \nonumber
\end{eqnarray}

Let $[\cdot\, , \cdot]_p$ be the \SIP\ of $L_p (I)$ given by
\begin{displaymath} 
                   [f,g]_p =  \int_I f \bigg ( \frac{\overline{g} \vert g \vert^{p-2}}{\Vert g \Vert_p^{p-2}}\bigg ) dm .
\end{displaymath}
Then, we can write the ``{band-limited}'' \UU\ \BSS\ and their duals as:
\begin{eqnarray}
                \BB_p   &:=& \{ f= [\widehat{f},\Phi (\cdot )]_p \in C(\RR)  : \widehat{f}\in L_p(I) \}  \nonumber\\
                \BB_p^* &:=& \{ h= [\Phi (\cdot ), \widehat{h}\,]_p \in C(\RR) : \widehat{h}\in L_p(I)  \}. \nonumber
\end{eqnarray}
By norming them with $\Vert f \Vert_{\BB_p}= \Vert \widehat{f}\, \Vert_{L_p (I)}$ and $\Vert h \Vert_{\BB_p^*}= \Vert \widehat{h}\, \Vert_{L_p (I)}$ respectively, 
indeed, we obtain two \UU\ \BSS. The duality between them can be written as
\begin{displaymath}
                   (f,h)_{\BB_p} = (\widehat{f}\,,(\widehat{h}\,)^*)_p = [\widehat{f}\,,\widehat{h}\,]_p
                                   \qquad f\in\BB_p,\ h\in\BB_p^*,
\end{displaymath} 
In these terms, $\BB_p$ is a \SIPRKBS\ on $\RR$ with the \SIP\ given by
\begin{equation}
                [f,g]_{\BB_p} = \big [[\widehat{f},\Phi (\cdot )]_p , [\widehat{g},\Phi (\cdot )]_p \big ]_{\BB_p} = [\widehat{f}\,, \widehat{g}\,]_p
                                \qquad f,g\in\BB_p\nonumber
\end{equation}
and the \SIPRK\ $G$ has the form
\begin{align*}
              G(\omega, t ) := [\Phi (\omega) , \Phi (t ) ]_p = 
              \int_I  e^{-2\pi i\omega x} e^{2\pi it x } dx =
              \frac{\sin \pi (t - \omega )}{\pi (t - \omega )} = sinc (t - \omega ) 
                          \qquad t,\omega\in \RR  ,
\end{align*}
The reproducing property is satisfied, since
\begin{align*}
              [ f, G(t,\cdot ) ]_{\BB_p} =
              \big[ [\widehat{f}\,,\Phi (\cdot)]_p,[\Phi (t) , \Phi (\cdot) ]_p \big ]_{\BB_p} =
              [\widehat{f}\,,\Phi (t )]_p =
              \int_{\RR} \widehat{f}(\omega ) e^{2\pi it\omega} d\omega =
              f(t)
\end{align*}
Consequently, we have $(\BB_p, [\cdot\,,\cdot ]_{\BB_p}, G,L_p (I), \Phi)$ a \SIPRKBS\ on $\RR$. 

Now, we consider the sequence $\{G_j (\cdot)\}_{j\in\ZZ}\in\BB_p$, where 
\begin{equation}
                G_j(t)=G(j,t)=sinc(t-j)\qquad j\in\ZZ,\qquad t\in\RR ,
\nonumber
\end{equation}
being the integers ordered by $\ZZ=\{0,-1,1,\cdots \}$. That $\{G_j (t)\}_{j\in\ZZ}\in\ell_q(\ZZ)\ \,\forall\,t\in\RR$ is due to
\begin{equation}
                \bigg (\sum_{j\in\ZZ}\vert G_j(t) \vert^q \bigg )^{\frac{1}{q}}\leq 
                \frac{1}{\pi}\bigg (\sum_{ j\in\ZZ}\frac{1}{\vert t-j \vert^q} \bigg )^{\frac{1}{q}} <
                +\infty .
\nonumber
\end{equation}
Of course, when $t\in\ZZ$, $G_t (t)=1$. The last calculation shows that $\{G_j (t)\}_{j\in\ZZ}\in\ell_q(\ZZ)\ \,\forall\,t\in\RR$ (and $\forall\,1<q<+\infty$ in fact),  
therefore $\{G_j (t)\}^*_{j\in\ZZ}\in\ell_p(\ZZ)\ \,\forall\,t\in\RR$. 

By calling $\phi:\RR\rightarrow\ell_q(\ZZ)$ to the map $t\mapsto \{G_j (t)\}_{j\in\ZZ}$ and $\phi^*:\RR\rightarrow\ell_p(\ZZ)$ to the map 
$t\mapsto \{G_j (t)\}^*_{j\in\ZZ}$, follows inmediately that
\begin{eqnarray}
                \overline{span}\big \{ \{ G_j (t)   \}_{j\in\ZZ}:t\in\RR\big\} & = & \ell_q(\ZZ)  \nonumber \\
                \overline{span}\big \{ \{ G_j (t) \}^*_{j\in\ZZ}:t\in\RR\big\} & = & \ell_p(\ZZ), \nonumber
\end{eqnarray}
since  $G_j(k)=\delta_{j,k}\ \,\forall\,j,k\in\ZZ$, thus the canonical unconditional basis for $\ell_q(\ZZ)$ (whenever 
$1<q<+\infty$) is contained in both sets. At this point, we already can define the \SIPRKBS\ on 
$\RR$ $(\BB_{samp} , [ \cdot ,\cdot ]_{samp}, G_{samp}, \ell_q(\ZZ), \phi )$, which \SIPRK\ $G_{samp}$ is given by 
\begin{equation}
                G_{samp}(s,t)=\sum_{j\in\ZZ}(G_j(s))^*G_j(t)= 
                \frac{1}{\Vert \{G_j(s)\}_j\Vert_{\ell_p(\ZZ)}^{p-2}}\sum_{j\in\ZZ}G_j(t)\overline{G_j(s)}\vert G_j(s)\vert^{p-2}.
\nonumber
\end{equation}
Also it holds true that $\BB_{samp}\subset \BB_p$ (by the last two completeness conditions).

We are going to see that $\{G_j (\cdot)\}_{j\in\ZZ}$ satisfies the hypotheses of the ``{Converse Sampling Theorem}'' so, we choose the sequences 
$\{t_j:=j\}_{j\in\ZZ}\subset\RR$ and $\{a_j:\equiv 1\}_{j\in\ZZ}\subset \CC\setminus \{0\}$. 

In the first place, we need to show that the sequence $\{f(j)\}_{j\in\ZZ}$ belongs to $\ell_p(\ZZ)$ for all $f\in\BB_p$. If $f\in\BB_p$ and $j\in\ZZ$, then
\begin{equation}
                f(j) = [f,G_j(\cdot )]_{\BB_p} = 
                [\widehat{f}\, , \Phi (t ) ]_p = 
                \int_{I}\widehat{f}\,(\omega)e^{2\pi ij\omega }d\omega =
                \int_{I}\widehat{f}\,(\omega )\widehat{G_j}\,(\omega )d\omega 
\nonumber
\end{equation}
so, for $j\neq 0$:
\begin{align*}
              \big\vert f(j) \big\vert^p  & =
              \bigg\vert  \int_{I}\widehat{f}\,(\omega)\widehat{G_j}\,(\omega )d\omega\bigg\vert^p \\ &=
              \bigg\vert  \int_{I}\widehat{f}\,(\omega)e^{2\pi ij\omega }d\omega\bigg\vert^p  \\ &=
              \bigg\vert \bigg [ \frac{\widehat{f}\,(\omega)e^{2\pi ij\omega }}{2\pi i j}\bigg ]^{\frac{1}{2}}_{-\frac{1}{2}} - 
                         \frac{1}{2\pi i j}\int_{I}\widehat{f}\,'(\omega)e^{2\pi ij\omega }d\omega\bigg\vert^p \\ &\leq
              \frac{C(p)\Vert\widehat{f}\,'\Vert^p_p}{\vert j\vert^p}.          
\end{align*}
while for $j=0$ we have $\big\vert f(j) \big\vert^p \leq \Vert\widehat{f}\,\Vert^p_p$. Therefore, taking the $\ell_p$-norm results:
\begin{displaymath}
                   \big \Vert \{ f(j) \}_{j\in\ZZ} \big \Vert_{\ell_p(\ZZ)}  \leq 
                   C(p,\widehat{f}\,,\widehat{f}\,') \bigg ( 1+2\sum_{j\in\NN} \frac{1}{j^p}\bigg )^{\frac{1}{p}} < +\infty .
\end{displaymath}

In the second place, we want the sampling representation $f(t)=\sum_{j\in\ZZ}f(j)G_j(t),\ t\in\RR$ for all $f\in\BB_p$, being the series pointwise 
convergent at least. The series in fact is absolutely convergent since $\{{f(j)}\}_{j\in\ZZ}\in\ell_p(\ZZ)$ and $\{G_j (t)\}_{j\in\ZZ}\in\ell_q(\ZZ)$ for all 
$f\in\BB_p,\ t\in\RR$. In this way only remains to check the pointwise convergence to $f(t)$. If $f\in\BB_p$, we have the following representation:
\begin{align*}
              f(t)=[ f, G(t,\cdot ) ]_{\BB_p} = \int_{I} \widehat{f}\,(\omega ) \widehat{G_t}\,(\omega )d\omega \qquad t\in\RR ,
\end{align*}
We consider the sequence of functions $\{ f_N \}_{N\in\NN_0}$ in $\BB_p$ given by  
\begin{equation}
                f_N(t)= \sum_{\vert j\vert\leq N} f(j)G_j(t)\qquad N\in\NN_0,\ t\in\RR.
\nonumber
\end{equation}
By fixing $f\in\BB_p$ and $t\in\RR$, for $N\in\NN_0$ we have:
\begin{align*}
              \vert f(t)-f_N(t)\vert &= 
              \bigg\vert\int_{I}\widehat{f}\,(\omega )\widehat{G_t}\,(\omega )d\omega - 
                                  \sum_{\vert j\vert\leq N}\bigg ( \int_{I}\widehat{f}\,(\omega )\widehat{G_j}\,(\omega )d\omega \bigg )G_j(t) \bigg\vert \\&\leq
              \int_{I} \vert\widehat{f}\,(\omega )\vert\bigg\vert \widehat{G_t}\,(\omega ) - \sum_{\vert j\vert\leq N}\widehat{G_j}\,(\omega )G_j(t)\bigg\vert d\omega \\&=
              \int_{I} \vert\widehat{f}\,(\omega )\vert\bigg\vert \sum_{ j\in\ZZ}G_t(j)\widehat{G_j}\,(\omega )
                                - \sum_{\vert j\vert\leq N}\widehat{G_j}\,(\omega )G_j(t)\bigg\vert d\omega \\ &= 
              \int_{I} \vert\widehat{f}\,(\omega )\vert\bigg\vert \sum_{\vert j\vert > N}\widehat{G_j}\,(\omega )G_j(t)\bigg\vert d\omega \\&\leq
              \Vert \widehat{f}\,\Vert_{L_p (I)} \bigg ( \int_I \bigg\vert \sum_{\vert j\vert > N}\widehat{G_j}\,(\omega )G_j(t)\bigg\vert^q d\omega \bigg )^{\frac{1}{q}} .
\end{align*}
Where in the third equality we used that $G_t(j)=G_j(t)\ \,\forall\,t\in\RR , j\in\ZZ$. Then, by Lebesgue's Dominated convergence Theorem, follows that $\{ f_N \}_{N\in\NN_0}$ 
converges pointwise to $f$ in $\RR$.  
\end{exam}  
%

\end{document}